# An asymptotically normal test for the selective neutrality hypothesis[*]

**Aluísio Pinheiro[1], Hildete P. Pinheiro[1] and Samara Kiihl[1]**

*Universidade Estadual de Campinas*

**Abstract:** An important parameter in the study of population evolution is $\theta = 4N\nu$, where $N$ is the effective population size and $\nu$ is the rate of mutation per locus per generation. Therefore, $\theta$ represents the mean number of mutations per site per generation. There are many estimators of $\theta$, one of them being the mean number of pairwise nucleotide differences, which we call $\mathcal{T}_2$. Other estimators are $\mathcal{T}_1$, based on the number of segregating sites and $\mathcal{T}_3$, based on the number of singletons. The concept of selective neutrality can be interpreted as a differentiated nucleotide distribution for mutant sites when compared to the overall nucleotide distribution. Tajima (1989) has proposed the so-called Tajima's test of selective neutrality based on $\mathcal{T}_2 - \mathcal{T}_1$. Its complex empirical behavior (Kiihl, 2005) motivates us to propose a test statistic solely based on $\mathcal{T}_2$. We are thus able to prove asymptotic normality under different assumptions on the number of sequences and number of sites via $U$-statistics theory.

## 1. Introduction

A large number of metrics has been constructed to measure genetic distances. Examples of such metrics are the Gini-Simpson index of diversity (Gini [8], Simpson [27] and Sen [23]), Nei, Mahalanobis and Hamming distances (Rao [20, 21] and Chakraborty and Rao [2]). Much work has been done using these measures of genetic distances to test homogeneity of genetic data or genetic polymorphism (Pinheiro et al. [19], Sen [24] and Pinheiro et al. [15]).

The Gini-Simpson index can be used to build a sub-additive analysis of variance for categorical data (Pinheiro et al. [17]). Moreover, tests of homogeneity among groups of genomic sequences can be constructed (Pinheiro et al. [15, 18] and Pinheiro et al. [19]) based on Hamming distance inequalities (Sen [23]).

In the study of population evolution, one is interested in the population's degree of polymorphism. This can be expressed as the rate of mutation per locus per generation. Some tests for selective neutrality include Tajima [29], Fu [5], Fu and Li [7] and Fu [6]. A communality to these tests is their intrinsic dependence on the Jukes and Cantor type of mutation processes (Jukes-Cantor [12]). However, there is abundant evidence for on the inadequacy of assumptions such as equal rates of mutation for every site (Fitch and Margoliash [4], Yang [34] and Uzzell and Corbin [32]).

[*]Supported in part by CNPq Grants 474329/2004-6, 476781/2004-3 and Fapesp Grant 2003/10105-2.

[1]Departamento de Estatística UNICAMP, Caixa Postal 6065 CEP 13083-970, Campinas SP, Brazil, e-mail: pinheiro@ime.unicamp.br; hildete@ime.unicamp.br

Another characteristic of genomic data is the challenge posed to parametric models by the enormously large number of sites, $K$, with relatively smaller sampled sequences, $n$. These statistical problems and some solutions for related measures are discussed by Sen et al. [26], Sen [25] and Pinheiro et al. [15].

In the study of population evolution, an important parameter of interest is $\theta = 4N\nu$, where $N$ is the effective population size and $\nu$ is the rate of mutation per locus per generation. Thus, $\theta$ represents the mean number of mutations per site per generation.

There are many estimators of $\theta$ in the literature (Hartl and Clark [9]). One of them is the *mean number of pairwise nucleotide differences*, which we will call $\mathcal{T}_2$. Other estimators of $\theta$ are $\mathcal{T}_1$, based on the *number of segregating sites* and $\mathcal{T}_3$, based on the *number of singletons*. The estimators $\mathcal{T}_1$, $\mathcal{T}_2$ and $\mathcal{T}_3$ are used in the literature to build test statistics for selective neutrality (Tajima [29] and Fu and Li [7]).

The concept of selective neutrality can be interpreted as a differentiated nucleotide distribution for mutant sites when compared to the overall nucleotide distribution. A statistic has been proposed, based on $\mathcal{T}_2 - \mathcal{T}_1$, for testing such a hypothesis (Tajima [29]). This, however, has shortcomings, due to the complex behavior of $\mathcal{T}_2 - \mathcal{T}_1$ (Kiihl [13]). We propose a test statistic solely based on $\mathcal{T}_2$. For this statistic, we are able to prove asymptotic normality under different asymptotic conditions on the number of sequences and sites and also derive Berry-Esséen rates of convergence.

The text goes as follows. Section 2 reviews the biological motivation and the tests for selective neutrality available in the literature. In Section 3 we propose a test solely based on the nucleotide frequencies. Its asymptotic behavior is studied in Section 4. In Section 5, we illustrate the performances of the proposed test and Tajima's procedure in two genetic data sets.

## 2. Tajima's Test of Selective Neutrality

In order to better understand the differences between these estimators, let us suppose that we have a sample of 5 DNA sequences, from which 500 sites were sequenced. Table 1 presents only the 16 polymorphic or *segregating sites*, i.e., those sites in which we find nucleotide differences. The other 484 sites which do not present differences are called *non-segregating sites*. Among the polymorphic sites, sites 3, 5, 6, 7, 8, 11, 12, 13 and 15 are *singletons*, since they present only one nucleotide different from the others.

Let $S$ be the number of *segregating* sites, $\mathcal{T}_2$ be the mean number of *pairwise differences* and $S^\star$ be the number of *singletons*. In the example given in Table 1,

TABLE 1
*Polymorphic sites in a sample of five genes*

|   | 1 | 2 | 3 | 4 | 5 | 6 | 7 | 8 | 9 | 10 | 11 | 12 | 13 | 14 | 15 | 16 |
|---|---|---|---|---|---|---|---|---|---|----|----|----|----|----|----|----|
| a | T | C | T | A | C | C | T | C | C | T | C | G | G | T | T | A |
| b | T | C | C | T | A | C | C | T | C | C | T | G | G | T | T | T |
| c | C | T | C | C | C | C | C | T | C | T | T | T | G | C | T | A |
| d | C | T | C | C | C | C | C | T | T | C | T | G | A | C | T | T |
| e | C | T | C | C | C | T | C | T | T | T | T | G | G | C | C | A |
| * | (6) | (6) | (4) | (7) | (4) | (4) | (4) | (4) | (6) | (6) | (4) | (4) | (4) | (6) | (4) | (6) |

\* : number of pairwise differences for each site.
There are 16 segregating sites.
Singletons are presented in boldface.

$S = 16$,

$$\mathcal{T}_2 = \frac{6+6+4+\cdots+6+4+6}{16} = 4.94$$

and $S^\star = 8$. The main difference between $S$ and $\mathcal{T}_2$ is the effect of selection.

The motivation and interpretation of the test of neutral mutation is that nucleotide polymorphism (segregating sites) and nucleotide diversity (pairwise nucleotide differences) differ primarily because the segregating sites are *indifferent* to the relative frequencies of the polymorphic nucleotides.

Mutant nucleotides are maintained in the population with low frequency. As the number of segregating sites ignores the frequency of mutant nucleotides, its value can be strongly affected by their existence, even if they occur with low frequency. On the other hand, the existence of mutant nucleotides with low frequency does not affect the mean number of pairwise differences, since in this case the frequency of mutations is considered. In other words, if some of the observed mutations have selective effects, the estimator $\mathcal{T}_1$ of $\theta$, based on $S$, cannot be the same as $\mathcal{T}_2$. Major discrepancies occur when:

- The relative frequencies of polymorphic variants are almost identical (nearly equal). This pattern increases the proportion of nucleotide pairwise differences, hence $\mathcal{T}_2 - \mathcal{T}_1$ is positive. This suggests either some type of balancing selection, in which heterozygous genotypes are favored, or some type of diversifying selection, in which genotypes carrying the less common alleles are favored;
- The relative frequencies of the polymorphic variants are too unequal, with an excess of the most common type and a deficiency of the less common types. This pattern results in a decrease in the proportion of pairwise differences, so $\mathcal{T}_2 - \mathcal{T}_1$ is negative. Typical reasons for excessively unequal frequencies can be selection against genotypes carrying the less frequent alleles, recent population bottleneck eliminating less frequent alleles, and insufficient time since the occurrence of the bottleneck to restore the equilibrium between mutation and random genetic drift.

The most common and well-known test of selective neutrality in the literature is Tajima's D test (Tajima [29]). This test uses the following statistic to test the hypothesis of neutral mutation, which is known as Tajima's $D$ statistic:

$$(2.1) \qquad D = \frac{D_1}{\sqrt{\mathrm{Var}(D_1)}},$$

where $D_1 = \mathcal{T}_2 - \mathcal{T}_1$, for which

$$(2.2) \qquad \mathrm{Var}(D_1) \approx \frac{n+1}{3(n-1)}\theta + \left[ \frac{2(n^2+n+3)}{9n(n-1)} - \frac{n+2}{a_n\, n} + \frac{b_n}{a_n^2} \right]\theta^2,$$

with $a_n = \sum_{j=1}^{n-1}(1/j)$ and $b_n = \sum_{j=1}^{n-1}(1/j^2)$.

Note that $D_1$ is not a $U$-statistic, i.e. $E(\mathcal{T}_1)$ is not an estimable parameter in Hoeffding's sense, since $\mathcal{T}_1$ estimates it by using the sample as a whole. Heuristic motivation is provided by Tajima [29] to use Beta distribution tables for the asymptotic behavior of $D$. In Kiihl [13] the theoretical characteristics of $D$ are carefully studied. It is shown that $D$'s asymptotic distribution is infinitely divisible but strong theoretical evidence is provided against asymptotic normality. In view of the difficulties in dealing with the theoretical distribution of $D$, the use of resampling methods, such as bootstrap, are recommended to generate its empirical distribution and to compute the p-value of the test.

### 3. The test of selective neutrality mutations based on nucleotide frequencies

Based on the interpretation of the vaguely defined alternative hypothesis discussed in Section 2, we can think of a hypothesis driven by differences on the nucleotide frequencies for segregating and non-segregating sites.

Let $\Pi_{cl}$ be the probability of having category (nucleotide) $c$ at site $l$, with $c = 1, \ldots, 4$; $l = 1, \ldots, K$ and $K$ is the number of segregating sites. Let $\Pi_c$ be the overall probability of having category (nucleotide) $c$ in non-segregating sites.

Now, let $\mathbf{X}_i = (X_{i1}, \ldots, X_{iK})'$ and $\mathbf{X}_j = (X_{j1}, \ldots, X_{jK})'$ be random vectors representing DNA sequences $i$ and $j$. So, $X_{il}$ can assume values in the set $\{A, C, T, G\}$, where $A$ represents Adenine, $C$, Cytosine, $T$, Thymine and $G$, Guanine. As in Pinheiro et al. [15] we write

$$(3.1) \qquad D_{ij} = \frac{1}{K} \sum_{l=1}^{K} \mathbb{I}(X_{il} \neq X_{jl})$$

the proportion of sites where $\mathbf{X}_i$ and $\mathbf{X}_j$ differ; here $\mathbb{I}(A)$ stands for the indicator function of set $A$. Consider

$$
\begin{aligned}
(3.2) \qquad \mathcal{H}_K &= \frac{1}{K} \sum_{l=1}^{K} P(X_{il} \neq X_{jl}) \\
&= \frac{1}{K} \sum_{l=1}^{K} \sum_{c=1}^{4} \Pi_{cl}(1 - \Pi_{cl}) = 1 - \frac{1}{K} \sum_{l=1}^{K} \sum_{c=1}^{4} \Pi_{cl}^2.
\end{aligned}
$$

A natural, unbiased and optimal nonparametric estimator of $\mathcal{H}_K$ is $\mathcal{T}_2$, a $U$-statistic (Hoeffding [10]) of degree 2, given by

$$(3.3) \qquad \mathcal{T}_2 = \frac{1}{K} \binom{n}{2}^{-1} \sum_{l=1}^{K} \sum_{1 \leq i < j \leq n} \mathbb{I}(X_{il} \neq X_{jl}).$$

We can write, under the null hypothesis of neutral mutation,

$$(3.4) \qquad H_0 : \mathrm{E}(\mathcal{T}_2) = 1 - \frac{1}{|\mathcal{N}|} \sum_{t \in \mathcal{N}} \sum_{c=1}^{4} \Pi_{ct}^2 \equiv \theta_0,$$

where $\mathcal{N}$ is the set of non-segregating sites and $|\mathcal{N}|$ is its cardinality.

Using the Hoeffding-decomposition (Hoeffding [10]), we have

$$(3.5) \qquad \mathcal{T}_2 = \mathcal{H}_K + 2H_n^{(1)} + H_n^{(2)},$$

where $H_n^{(1)} = n^{-1} \sum_{i=1}^{n} E(D_{ij} \mid X_j) - \mathcal{H}_K$, and $H_n^{(2)} = \mathcal{T}_2 - 2H_n^{(1)} + \mathcal{H}_K$ is the degenerate component of order 2.

If $\sigma_1^2$, the variance of $E(D_{ij} \mid X_j)$, is positive, $\mathcal{T}_2$ is a nondegenerate $U$-statistic of degree 2, and by Hoeffding [10],

$$(3.6) \qquad n^{1/2}(\mathcal{T}_2 - \mathcal{H}_K) \overset{D}{\longrightarrow} \mathrm{N}(0, 4\sigma_1^2).$$

So, in our case, if we assume the same conditions as in Tajima [29], that is, independence among sites and sequences, we get asymptotic normality for the test statistic $\mathcal{T}_2$, under either $H_0$ or $H_1$. Moreover, since the kernel is bounded, usual fourth moment conditions easily hold and Berry-Esséen results are straight forward albeit tedious and cumbersome. More powerful tools are provided in Pinheiro et al. [16, 19] and Pinheiro et al. [15], which we can somewhat mimic to relax some of the initial conditions of Tajima [29].

## 4. Asymptotics of the neutral selectivity test

Motivated by the discussion in Section 3, we define

$$(4.1) \qquad \mathcal{T}_2 = \frac{1}{K\binom{n}{2}} \sum_{(i,j)} \sum_{k=1}^{K} \mathbb{I}(X_{ik} \neq X_{jk}).$$

$\mathcal{T}_2$ is a non-degenerate $U$-statistic of degree 2, for which

$$(4.2) \qquad \phi(\mathbf{X}_i, \mathbf{X}_j) = \frac{1}{K} \sum_{k=1}^{K} \mathbb{I}(X_{ik} \neq X_{jk}),$$

$$(4.3) \qquad \psi_1(\mathbf{X}_1) = 1 - \frac{1}{K} \sum_{k=1}^{K} \Pi_{X_{1k}k}$$

$$(4.4) \qquad h^{(1)}(\mathbf{X}_1) = \frac{1}{K} \sum_{k=1}^{K} \left[ \Pi_{X_{1k}k} - \sum_{c=1}^{C} \Pi_{ck}^2 \right]$$

$$h^{(2)}(\mathbf{X}_1, \mathbf{X}_2) = \frac{1}{K} \sum_{k=1}^{K} \mathbb{I}(X_{ik} \neq X_{jk}) - 1 + \frac{1}{K} \sum_{k=1}^{K} \Pi_{X_{1k}k}$$

$$(4.5) \qquad \qquad + \frac{1}{K} \sum_{k=1}^{K} \Pi_{X_{2k}k} - \frac{1}{K} \sum_{k=1}^{K} \sum_{c=1}^{C} \Pi_{ck}^2,$$

where $\Pi_{ck} = P(X_{ik} = c)$, for $k = 1, \ldots, K$ and $c = 1, \ldots, C$.

Note that

$$\sigma_1^2 = Eh^{(1)}(\mathbf{X}_1)^2$$

$$= \frac{1}{K^2} \left\{ \sum_{k=1}^{K} \sum_{c=1}^{C} \Pi_{ck}^3 + \sum_{k \neq l=1}^{K} \sum_{c,d=1}^{C} (\Pi_{ck,dl} - \Pi_{ck}\Pi_{dl}) \Pi_{ck}\Pi_{dl} \right.$$

$$(4.6) \qquad \left. - \sum_{k=1}^{K} \sum_{c,d=1}^{C} \Pi_{ck}^2 \Pi_{dk}^2 \right\}.$$

Suppose, for instance, that the $K$ sites are independently distributed. Then (4.6) reduces to

$$\sigma_1^2 = \frac{1}{K^2} \left\{ \sum_{k=1}^{K} \sum_{c=1}^{C} \Pi_{ck}^3 - \sum_{k=1}^{K} \sum_{c,d=1}^{C} \Pi_{ck}^2 \Pi_{dk}^2 \right\},$$

which will be zero if and only if $\Pi_{ck} = 1/C$, $c = 1, \ldots, C$ and $k = 1, \ldots, K$ or, if for each $k = 1, \ldots, K \; \exists c \in \{1, \ldots C\}$ such that $\Pi_{ck} = 1$. We can then assure that $\mathcal{T}_2$ is a non-degenerate U-statistic of degree 2 unless either all sites' distributions are degenerate or uniform, which can be generally classified as non-interesting cases for genetic (or otherwise) data. If, however, one is interested in such null hypotheses, we refer the reader to Pinheiro et al. [15] and Pinheiro et al. [19] for related issues in a somewhat different approach. For a more general setup in which both null and alternative hypotheses have associated generalized degenerate (quasi) U-statistics, we refer the reader to Pinheiro et al. [16].

We also know that

$$
\begin{aligned}
\mathcal{H}_K &\equiv E(\mathcal{T}_2) = K^{-1} \binom{n}{2}^{-1} \binom{n}{2} \sum_{k=1}^{K} \sum_{c=1}^{C} \Pi_{ck}(1 - \Pi_{ck}) \\
&= 1 - \frac{1}{K} \sum_{k=1}^{K} \sum_{c=1}^{C} \Pi_{ck}^2
\end{aligned}
$$

(4.7)

and its H-decomposition is given by

$$
\mathcal{T}_2 = \mathcal{H}_K + \frac{2}{n} \sum_{i=1}^{n} h^{(1)}(\mathbf{X}_i) + \binom{n}{2}^{-1} \sum_{i<j} h^{(2)}(\mathbf{X}_i, \mathbf{X}_j).
$$

We can then write

$$
\sqrt{K}n \frac{\mathcal{T}_2 - \mathcal{H}_K}{2\sigma_1} = \sum_{i=1}^{n} \frac{\sqrt{K}}{\sqrt{n}\sigma_1} h^{(1)}(\mathbf{X}_i) + o_p(1),
$$

i.e., the asymptotic behavior of $\sqrt{K}n(\mathcal{T}_2 - \mathcal{H}_K)$ depends only on the asymptotic behavior of the sum

$$
Z_n = \sum_{i=1}^{n} \frac{1}{\sqrt{nK}\sigma_1} \sum_{k=1}^{K} \left[ \Pi_{X_{i_k}k} - \sum_{c=1}^{C} \Pi_{ck}^2 \right].
$$

If we write $Z_n = \sum_{i=1}^{n} Y_{ni}$, where

$$
Y_{ni} = \frac{1}{\sqrt{nK}\sigma_1} \sum_{k=1}^{K} \left[ \Pi_{X_{i_k}k} - \sum_{c=1}^{C} \Pi_{ck}^2 \right],
$$

we need to show the CLT for the array $\{Y_{ni}, n \geq 1, 1 \leq i \leq n\}$, for $k \equiv k(n)$.

If $K$ is finite, $\mathcal{H}_K \equiv \mathcal{H}$ and CLT's for r.v.'s will suffice. Since $Var(\mathcal{T}_2) > 0$, Hoeffding [10]'s result can be applied and

(4.8)                    $\sqrt{n}K \dfrac{\mathcal{T}_2 - \mathcal{H}}{2\sigma_1} \xrightarrow{D} N(0,1)$ as $n \to \infty$.

If, however, $K$ varies with $n$, note that

$$
\begin{aligned}
s_n^2 &= \sum_{i=1}^{n} E Y_{ni}^2 \\
&= n \frac{1}{nK\sigma_1^2} E \left\{ \sum_{k=1}^{K} \left[ \Pi_{X_{i_k}k} - \sum_{c=1}^{C} \Pi_{ck}^2 \right] \right\}^2.
\end{aligned}
$$

Under mixing conditions (such as $\alpha(k) = \gamma^{-k}$, $0 < \gamma < 1$), $E\left\{\sum_{k=1}^{K}\left[\Pi_{X_{i_k}k} - \sum_{c=1}^{C}\Pi_{ck}^2\right]\right\}^2 = O(K)$ as $K \to \infty$, and therefore $s_n^2 = O(1)$ as $n \to \infty$ (and either $K$ limited or $K \to \infty$). On the other hand,

$$
\begin{aligned}
E|Y_{ni}|^4 &= (\sigma_1^2 nK)^{-4/2} E\left|\sum_{k=1}^{K}\left[\Pi_{X_{i_k}k} - \sum_{c=1}^{C}\Pi_{ck}^2\right]\right|^4 \\
&= (\sigma_1^2 nK)^{-2} E\sum_{k=1}^{K}\left[\Pi_{X_{i_k}k} - \sum_{c=1}^{C}\Pi_{ck}^2\right]\sum_{l=1}^{K}\left[\Pi_{X_{i_l}l} - \sum_{c=1}^{C}\Pi_{cl}^2\right] \\
&\quad \times \sum_{m=1}^{K}\left[\Pi_{X_{i_m}m} - \sum_{c=1}^{C}\Pi_{cm}^2\right]\sum_{p=1}^{K}\left[\Pi_{X_{i_p}p} - \sum_{c=1}^{C}\Pi_{cp}^2\right] \\
&= (\sigma_1^2 nK)^{-2}\left\{\sum_{k=1}^{K}\sum_{d=1}^{C}\Pi_{dk}\left(\Pi_{dk} - \sum_{c=1}^{C}\Pi_{ck}^2\right)^4\right. \\
&\quad + 6\sum_{k \neq l=1}^{K}\sum_{d,e=1}^{C}\Pi_{dk,el}\left(\Pi_{dk} - \sum_{c=1}^{C}\Pi_{ck}^2\right)^2\left(\Pi_{el} - \sum_{c=1}^{C}\Pi_{cl}^2\right)^2 \\
&\quad + 4\sum_{k \neq l \neq m=1}^{K}\sum_{d,e,f=1}^{C}\Pi_{dk,el,fm}\left(\Pi_{dk} - \sum_{c=1}^{C}\Pi_{ck}^2\right)^3\left(\Pi_{el} - \sum_{c=1}^{C}\Pi_{cl}^2\right) \\
&\quad + 3\sum_{k \neq l=1}^{K}\sum_{d,e=1}^{C}\Pi_{dk,el}\left(\Pi_{dk} - \sum_{c=1}^{C}\Pi_{ck}^2\right)^2\left(\Pi_{el} - \sum_{c=1}^{C}\Pi_{cl}^2\right)^2 \\
&\quad + 6\sum_{k \neq l \neq m=1}^{K}\sum_{d,e,f=1}^{C}\Pi_{dk,el,fm}\left(\Pi_{dk} - \sum_{c=1}^{C}\Pi_{ck}^2\right)^2\left(\Pi_{el} - \sum_{c=1}^{C}\Pi_{cl}^2\right) \\
&\quad \times \left(\Pi_{fm} - \sum_{c=1}^{C}\Pi_{cm}^2\right) \\
&\quad + \sum_{k \neq l \neq m \neq p=1}^{K}\sum_{d,e,f,g=1}^{C}\Pi_{dk,el,fm,gp}\left(\Pi_{dk} - \sum_{c=1}^{C}\Pi_{ck}^2\right)\left(\Pi_{el} - \sum_{c=1}^{C}\Pi_{cl}^2\right) \\
&\quad \left. \times \left(\Pi_{fm} - \sum_{c=1}^{C}\Pi_{cm}^2\right)\left(\Pi_{gp} - \sum_{c=1}^{C}\Pi_{cp}^2\right)\right\}.
\end{aligned}
$$

Using the same mixing rate, $\alpha(k) = \gamma^{-k}$, $E|Y_{ni}|^4 = O(n^{-2})$ as $n \to \infty$ (and either $K$ limited or $K \to \infty$). Therefore,

$$
\left(\sum_{i=1}^{n}E|Y_{ni}|^4\right)^2 = O(n^{-2}) = o((s_n^2)^4). \tag{4.9}
$$

By (4.9), Liapounov's CLT necessary conditions are attained and, therefore,

$$
\sqrt{nK}\frac{\mathcal{T}_2 - \mathcal{H}_K}{2\sigma_1} \xrightarrow{D} N(0,1) \tag{4.10}
$$

as $n \to \infty$ (and either $K$ limited or $K \to \infty$),

following Utev [31].

We should note that the convergence given by (4.10) is true also if $K \to \infty$ but $n$ is limited, simply using the CLT for mixing r.v.'s (Withers [33]). Therefore, we have asymptotic normality for $n \to \infty$ and/or $K \to \infty$. The only (sufficient) condition for that is mixing along the sequences if $K$ is large. Anyhow, such a hypothesis is much less restrictive than the sequencewise independence, equal rates of mutation for all sites, *infinite* number of sites, and Poisson distribution for the number of mutant sites taken by Tajima [29] and yet not sufficient for the asymptotic normality of $D$.

We can also assess the rate of convergence of the $U$-statistics' CLT by the appropiate generalization of Berry-Esséen's original results. For instance, if $n \to \infty$ and $K$ finite, Korolyuk and Borovskikh [14] proves that

$$
(4.11) \quad
\begin{aligned}
&\left| P\left( \sqrt{n} \frac{\mathcal{T}_2 - \mathcal{H}_K}{\sqrt{Var\mathcal{T}_2}} \leq x \right) - \Phi(x) \right| \\
&\leq C \left( \sigma_1^{-3} E |h^{(1)}(\mathbf{X}_1)|^3 + \sigma_1^{-5/3} E |h^{(2)}(\mathbf{X}_1, \mathbf{X}_2)|^{5/3} \right) n^{-1/2}.
\end{aligned}
$$

In our case, $\sigma_1^2$ is given by (4.6), $E|h^{(1)}(\mathbf{X}_1)|^3 \leq \left( E h^{(1)}(\mathbf{X}_1)^4 \right)^{3/4}$, which is given by (4.12), and $E|h^{(2)}(\mathbf{X}_1, \mathbf{X}_2)|^{5/3} \leq \left( E h^{(2)}(\mathbf{X}_1, \mathbf{X}_2)^2 \right)^{5/6}$, which is given by (4.13). We can write

$$
\begin{aligned}
&E h^{(1)}(\mathbf{X}_1)^4 \\
&= \frac{1}{K^4} \left\{ \sum_{k=1}^{K} \sum_{d=1}^{C} \Pi_{dk} \left( \Pi_{dk} - \sum_{c=1}^{C} \Pi_{ck}^2 \right)^4 \right. \\
&\quad + 3 \sum_{k \neq l=1}^{K} \sum_{d,e=1}^{C} \Pi_{dk,el} \left( \Pi_{dk} - \sum_{c=1}^{C} \Pi_{ck}^2 \right)^2 \\
&\quad \times \left( \Pi_{el} - \sum_{c=1}^{C} \Pi_{cl}^2 \right)^2 + 4 \sum_{k \neq l=m=1}^{K} \sum_{d,e,f=1}^{C} \Pi_{dk,el,fm} \left( \Pi_{dk} - \sum_{c=1}^{C} \Pi_{ck}^2 \right)^3 \\
&\quad \times \left( \Pi_{el} - \sum_{c=1}^{C} \Pi_{cl}^2 \right) \\
&\quad + 3 \sum_{k \neq l=1}^{K} \sum_{d,e=1}^{C} \Pi_{dk,el} \left( \Pi_{dk} - \sum_{c=1}^{C} \Pi_{ck}^2 \right)^3 \left( \Pi_{el} - \sum_{c=1}^{C} \Pi_{cl}^2 \right) \\
&\quad + 6 \sum_{k \neq l \neq m=1}^{K} \sum_{d,e,f=1}^{C} \Pi_{dk,el,fm} \left( \Pi_{dk} - \sum_{c=1}^{C} \Pi_{ck}^2 \right)^2 \left( \Pi_{el} - \sum_{c=1}^{C} \Pi_{cl}^2 \right) \\
&\quad \times \left( \Pi_{fm} - \sum_{c=1}^{C} \Pi_{cm}^2 \right) + \sum_{k \neq l \neq m \neq p=1}^{K} \sum_{d,e,f,g=1}^{C} \Pi_{dk,el,fm,gp} \left( \Pi_{dk} - \sum_{c=1}^{C} \Pi_{ck}^2 \right) \\
(4.12) \quad &\quad \left. \times \left( \Pi_{el} - \sum_{c=1}^{C} \Pi_{cl}^2 \right) \left( \Pi_{fm} - \sum_{c=1}^{C} \Pi_{cm}^2 \right) \left( \Pi_{gp} - \sum_{c=1}^{C} \Pi_{cp}^2 \right) \right\}
\end{aligned}
$$

and

$$Eh^{(2)}(X_1, X_2)^2$$

$$= \frac{1}{K^2} \left\{ K\mathcal{H}_K - \sum_{k=1}^{K} \sum_{c,d=1}^{C} \Pi_{ck}\Pi_{dk}(1-\Pi_{ck})(1-\Pi_{dk}) \right.$$

$$+ 2\sum_{k=1}^{K} \sum_{c,d=1}^{C} \Pi_{ck}^2\Pi_{dk}^2 - 2\sum_{k=1}^{K} \sum_{c=1}^{C} \Pi_{ck}^3$$

$$\left. + \sum_{k \neq l=1}^{K} \sum_{c,d=1}^{C} (\Pi_{ck,dl} - \Pi_{ck}\Pi_{dl})[1-\Pi_{ck} - \Pi_{dk} + \Pi_{ck,dl} - \Pi_{ck}\Pi_{dl}] \right\}.$$

(4.13)

For independent sites, (4.12) reduces to

$$Eh^{(1)}(X_1)^4 = \frac{1}{K^4} \left\{ \sum_{k=1}^{K} \sum_{d=1}^{C} \Pi_{dk} \left( \Pi_{dk} - \sum_{c=1}^{C} \Pi_{ck}^2 \right)^4 + \right.$$

$$+ 3\left[ \left\{ \sum_{k=1}^{K} \sum_{d=1}^{C} \Pi_{dk} \left( \Pi_{dk} - \sum_{c=1}^{C} \Pi_{ck}^2 \right)^2 \right\}^2 \right.$$

$$\left. \left. - \sum_{k=1}^{K} \left\{ \sum_{d=1}^{C} \Pi_{dk} \left( \Pi_{dk} - \sum_{c=1}^{C} \Pi_{ck}^2 \right)^2 \right\}^2 \right] \right\}$$

(4.14)

and (4.13) becomes

$$Eh^{(2)}(X_1, X_2)^2$$

$$= \frac{1}{K^2} \left\{ K\mathcal{H}_K - \sum_{k=1}^{K} \sum_{c,d=1}^{C} \Pi_{ck}\Pi_{dk}(1-\Pi_{ck})(1-\Pi_{dk}) \right.$$

$$\left. + 2\sum_{k=1}^{K} \sum_{c,d=1}^{C} \Pi_{ck}^2\Pi_{dk}^2 - 2\sum_{k=1}^{K} \sum_{c=1}^{C} \Pi_{ck}^3 \right\}.$$

(4.15)

Likewise, if $K \to \infty$ and finite $n$, under mixing conditions, we get

$$\left| P\left( \sqrt{K}\frac{\mathcal{T}_2 - \mathcal{H}_K}{\sqrt{Var\mathcal{T}_2}} \leq x \right) - \Phi(x) \right| \leq C\eta(K),$$

(4.16)

where $\eta(K)$ will be $K^{-1/2}\log K$ or slower than $K^{-1/2}$ depending on the mixing rates being exponential or polynomial on $K$, respectively, since the random variables are all bounded (Tihomirov [30]). Under site independence, the rate will be $K^{-1/2}$ (Feller [3]).

$H_0$, defined by (3.4), can be tested by

$$T_n = \sqrt{nK}\frac{\mathcal{T}_2 - \theta_0}{\sqrt{Var\mathcal{T}_2}}.$$

(4.17)

Note that under $H_0$, $T_n$ is asymptotically $N(0, 1)$ while, under $H_1$, $T_n = Z_n + \sqrt{nK}(\mathcal{H}_K - \theta_0)/\sqrt{Var\mathcal{T}_2}$, where $Z_n$ is asymptotically $N(0, 1)$ and $T_n - Z_n = O(nK)$. Therefore, when $n \to \infty$ and/or $K \to \infty$, the test defined by rejecting $H_0$ when $T_n > q_\alpha$ has an asymptotic size $\alpha$ and an asymptotic power one for all alternatives for which $\sqrt{nK}(\mathcal{H}_K - \theta_0)/\sqrt{Var\mathcal{T}_2} \to \infty$: that will happen whenever $K$ is finite and $\mathcal{H}_K - \theta_0 > 0$ or when $K \to \infty$ and $0 < \mathcal{H}_K - \theta_0 = O(n^{-1}K^{-1}\zeta(n, K))$, where $\zeta(n, K) \to \infty$. We can also deal with Pitman alternatives, i.e, for which $0 < \mathcal{H}_K - \theta_0 = O(n^{-1}K^{-1})$. Analogous reasonings will work for the test defined by $T_n < -q_\alpha$ or $|T_n| > q_{\alpha/2}$.

## 5. Applications

We illustrate the test statistic defined by (4.1) in two data sets. The first is composed by sequences of *Hydromedusa maximiliani*, a neotropical freshwater turtle from the Atlantic Forest's rivers in southeast Brazil (Souza et al. [28]). The sample is composed by $n = 48$ sequences of $K = 262$bp on the b cytochrome mitochondrial DNA. The second data set is composed by $n = 12$ HIV sequences from a single infected patient (Holmes and Brown [11]). Each sequence has $K = 233$bp.

For the computation of the variance of $\mathcal{T}_2$, we employed its jackknife estimator (Sen [22] and Arvesen [1]):

$$(5.1) \qquad \widehat{Var(\mathcal{T}_2)} = n^2(n-1)\binom{n-1}{m}^{-2} \sum_{c=0}^{m}(cn - m^2)S_c,$$

where $S_c = \sum \phi(\mathbf{X}_{i_1}, \mathbf{X}_{i_2})\phi(\mathbf{X}_{i_3}, \mathbf{X}_{i_4})$, for any resample $\{i_1, i_2, i_3, i_4\}$ from $\{1, \ldots, n\}$ such that there are $c$ coincident indices; $c = 0, \ldots, m$.

For the turtle data set, we find an asymptotic p-value of $2.021 \times 10^{-259}$ (test statistic equals $-34.41$). For the HIV data set, we find an asymptotic p-value of $3.131 \times 10^{-20}$ (teste statistic equals $-9.14$). Therefore, there is very strong statistical evidences for negative selection for both the turtle and the HIV data sets. On the other hand, if one uses Tajima's D statistic, one will decide on negative selection for the turtle data but on neutral selectivity for the HIV data with p-values 0.0397 and 0.7298, respectively (Kiihl [13]), even though for both data sets, the observed D statistic is negative.

In order to understand these test results, we should recall that the literature's accepted interpretation of neutral selectivity is that there is not much *change* on the sitewise nucleotide distribution for the mutant sites when compared to the non-mutant ones. Likewise, negative (positive) selectivity is interpreted as a significant bias towards a more concentrated nucleotide distribution (towards the uniform discrete distribution) for the segregating sites when compared to the non-segregating ones. We will take the observed frequencies as $\hat{\mathbf{\Pi}} = (\hat{\Pi}_A, \hat{\Pi}_C, \hat{\Pi}_G, \hat{\Pi}_T)'$. For the turtle data set, the non-segregating sites have observed frequencies given by $\hat{\mathbf{\Pi}}_{TUR,NS} = (0.3913, 0.2727, 0.2372, 0.0988)'$ while $\hat{\mathbf{\Pi}}_{TUR,S} = (0.8571, 0.1429, 0.0000, 0.0000)'$ for the segregating sites. For the HIV data set, the non-segregating sites have $\hat{\mathbf{\Pi}}_{HIV,NS} = (0.4550, 0.1327, 0.1706, 0.2417)'$ and we get $\hat{\mathbf{\Pi}}_{HIV,S} = (0.4091, 0.1364, 0.4545, 0.0000)'$ for the segregating sites. One should also note that $(1 - \hat{\mathbf{\Pi}}'_{TUR,S}\hat{\mathbf{\Pi}}_{TUR,S}) - (1 - \hat{\mathbf{\Pi}}'_{TUR,NS}\hat{\mathbf{\Pi}}_{TUR,NS}) = -0.4615$ and $(1 - \hat{\mathbf{\Pi}}'_{HIV,S} \times \hat{\mathbf{\Pi}}_{HIV,S}) - (1 - \hat{\mathbf{\Pi}}'_{HIV,NS}\hat{\mathbf{\Pi}}_{HIV,NS}) = -0.0804$, which clearly point to negative selectivity in both datasets.

To show the behavior of the test statistic under positive selectivity, we have artificially substituted the nucleotide frequencies for the seven segregating sites on the turtle data set by nucleotides closer to uniformly distributed ones. We have for the segregating sites an overall nucleotide frequency of .25, varying sitewise from 11/48 to 15/48 for each nucleotide class. The test value is then 5.7348 with a p-value of $9.7627 \times 10^{-9}$. We again reject the null hypothesis of selective neutrality, this time from a positive selectivity point of view.

One can also notice large changes in the observed frequencies between segregate and non-segregate sites in both examples, which agrees with the proposed test's conclusions. Tajima's test, however, leads to a different result for the HIV data. We should also point out that, since the jackknifed estimate of variance is positively biased (Sen [22] and Arvesen [1]), our results are conservative towards the null hypotheses, which strengthens even more the superior performance of the proposed test over Tajima's for both data sets.

## 6. Conclusions

We presented the inherent flaws associated with Tajima's test for selective neutrality, among them the vague definition of the null hypothesis (neutral selectivity), the (possibly) non-normal asymptotic distribution of $D$, the small number of presumed independent sites and the theoretically large number of sequences. We proposed a null hypothesis which formalizes the vague notions contained in Tajima's ideas for neutral selectivity. Moreover, due to $U$-statistics H-decomposition, we are able to provide our test with normal asymptotics, under a broader setup than those considered in the literature. The attained relaxations include: mixing positions (if $K$ is large) or any dependence setup (if $K$ is finite) instead of independently distributed positions; large $n$ and/or $K$ for normal asymptotics instead of large $K$ for motivation but only large $n$ for non-normal asymptotics. Resampling schemes which can be possibly cumbersome are easily circumvented by a direct formula for jackknifed $U$-statistics proposed by Sen [22] and Arvesen [1]. We illustrate the superior performance of the proposed test statistic to Tajima's in two data sets. In one of the data sets, we get a different conclusion (which is more reasonable when looking at other descriptive statistics). For the other data set, we come to the same conclusion, but the p-value is much smaller which is again biologically and statistically more reasonable. Summarizing, the proposed test uses all the advantages of $U$-statistics asymptotics, can be employed in a more general setup, and its application is quite simple due to jackknife variance estimation.

**Acknowledgments.** The authors would like to thank the editors and reviewers for their helpful and insightful suggestions which yield a clearer and richer text.